\documentclass[a4paper,11pt]{article}
\usepackage{amssymb,amsthm,amsmath,latexsym}

\newtheorem{thm}{Theorem}[section]
\newtheorem{pro}[thm]{Proposition}
\newtheorem{lem}[thm]{Lemma}
\newtheorem{cor}[thm]{Corollary}

\newtheorem*{thmA}{Theorem A}
\newtheorem*{thmB}{Theorem B}

\date{}

\title{Engel groups with an identity}

\author{\small{\textsc{Pavel Shumyatsky}\footnote{The first author was  supported by FAPDF/Brazil.}}\\
\small{Department of Mathematics, University of Brasilia}\\
\small{Brasilia-DF, 70910-900 Brazil}\\
\small{E-mail: pavel@unb.br}\\
[10pt]
\small{\textsc{Antonio Tortora}\footnote{The last two authors are members of {\em National Group for Algebraic and Geometric Structures, and their Applications} (GNSAGA--INdAM).} \;and \textsc{Maria Tota}\footnotemark[2]}\\
\small{Dipartimento di Matematica, Universit\`a di Salerno}\\
\small{Via Giovanni Paolo II, 132 - 84084 - Fisciano (SA), Italy}\\
\small{E-mail: antortora@unisa.it, mtota@unisa.it}}

\begin{document}
\maketitle

\begin{abstract} 
We give an affirmative answer to the question whether a residually finite Engel group satisfying an identity is locally nilpotent. More generally, for a residually finite group $G$ with an identity, we prove that the set of right Engel elements of $G$ is contained in the Hirsch-Plotkin radical of $G$. Given an arbitrary word $w$, we also show that the class of all groups $G$ in which the $w$-values are right $n$-Engel and $w(G)$ is locally nilpotent is a variety.\\

\end{abstract}

\section{Introduction}

Let $G$ be a group. An element $x\in G$ is called a right Engel element if for any $g\in G$ there exists a positive integer $n=n(x,g)$ such that $[x,_n g]=1$, where the commutator $[x,_n g]$ is defined inductively by the rules
$$[x,_1 g]=[x, g]=x^{-1}x^g\quad {\rm and,\, for}\; n\geq 2,\quad [x,_n g]=[[x,_{n-1} g],g].$$
Similarly, $x$ is a left Engel element if the variable $g$ appears on the left. 
If $n$ can be chosen independently of $g$, then $x$ is a right or left $n$-Engel element, respectively. 

A group $G$ is called an Engel group (or $n$-Engel, resp.) if its elements are both left and right Engel (or $n$-Engel, resp.). A theorem of Wilson states that a residually finite $n$-Engel group is locally nilpotent \cite{wi91}. On the other hand, the famous example of Golod shows that a residually finite Engel group is not in general locally nilpotent (see \cite{Rob}). 

Let $w=w(x_1,\ldots,x_d)$ be a nonempty word in the free group generated by $x_1,\ldots,x_d$. For a group $G$, we denote by $w(G)$ the verbal subgroup of $G$ corresponding to $w$, that is, the subgroup generated by the set $\{w(g_1, \dots ,g_d)\,|\,g_i\in G\}$ of all $w$-values
in $G$. Also, we say that $G$ satisfies the identity $w\equiv 1$ if $w(g_1,\ldots,g_d)=1$ for all $g_1,\ldots,g_d\in G$. The class of all groups satisfying the identity $w\equiv 1$ is called the variety determined by $w$. By a well-known theorem of Birkhoff, varieties are precisely classes of groups closed with respect to taking subgroups, quotients and Cartesian products of their members.

The main purpose of this paper is to prove the following theorem which confirms Conjecture 1.1 of \cite{BMTT}.

\begin{thmA}
Let $G$ be a residually finite group satisfying an identity. Then the set of right Engel elements of $G$ is contained in the Hirsch-Plotkin radical of $G$. In particular, if $G$ is an Engel group, then $G$ is locally nilpotent.	
\end{thmA}
	
Recall that the Hirsch-Plotkin radical of $G$ is the unique maximal normal locally nilpotent subgroup containing all normal locally nilpotent subgroups of $G$ (see \cite[12.1.3]{Rob}). The proof of Theorem A is based on Lie-theoretic techniques created by Zelmanov in his solution of the restricted Burnside problem. Other results in the same spirit were obtained in \cite{BSTT13,STT16} (see also \cite{STT15}) for left Engel elements.

An interesting application of Theorem A is the next result concerning varieties of groups with an Engel identity. 

\begin{thmB}\label{thmB} Let $n$ be a positive integer and $w$ an arbitrary word. Then the class of all groups in which the $w$-values are right $n$-Engel and $w(G)$ is locally nilpotent is a variety.
\end{thmB}

Theorem B may be compared to \cite[Theorem A]{STT16_2}, where the same authors proved that, for a {\em multilinear commutator} word $w$, the class of all groups $G$ in which the $w$-values are left $n$-Engel and $w(G)$ is locally nilpotent is a variety.

\section{About Lie algebras with an identity}

In this section we collect some definitions and results on Lie algebras satisfying an identity.

Let $L$ be a Lie algebra over a field and let $a,b_1,\dots,b_n$ be elements of $L$. We use the left normed convention for Lie brackets, that is,
$$[a,b_1,\dots,b_n]=[[\dots[[a,b_1],b_2],\dots],b_n];$$
an useful notation is 
$$[a,_n b_i]=[a,\underbrace{b_i,\dots,b_i}_{n \ times}].$$
An element $b\in L$ is called ad-nilpotent if there exists a positive integer $n$ such that $[a,_n b]=0$ for all $a\in L$. If $n$ is the least integer with this property, then $b$ is ad-nilpotent of index $n$. 
Following \cite{ze16}, we say that a subset $X$ of $L$ is a Lie set if $[a,b]\in X$ for any $a,b\in X$. We denote by $S\langle X\rangle$ the Lie set generated by $X$, namely the smallest Lie set containing $X$.

Let $F$ be the free Lie algebra over the same field as $L$ on the generators $x_1,\dots,x_m$. For a nonzero element $f=f(x_1,\dots,x_m)$ of $F$, the Lie algebra $L$ is said to satisfy the polynomial identity $f\equiv 0$ if $f(a_1,\dots,a_m)=0$ for all $a_1,\dots,a_m\in L$. 

As we will see in Section 3, Theorem A depends crucially on the following result of Zelmanov (\cite[Theorem 1.1]{ze16}; see also \cite{ze1,ze2}).

\begin{thm}\label{Zelmanov}
Let $L$ be a Lie algebra satisfying a polynomial identity and generated by elements $a_1,\ldots, a_m$. If every element $b\in S\langle a_1,\ldots,a_m\rangle$ is ad-nilpotent, then $L$ is nilpotent.
\end{thm}

Using \cite[Lemma 5]{MZ} we deduce the following corollary.

\begin{cor}\label{MZ}
For $p$ a prime, let $L$ be a Lie algebra over the field with $p$ elements satisfying a polynomial identity. Suppose that $L$ is generated by elements $a_1,\ldots, a_m$ such that  $[a_i,_{p^k} b]=0$, for some $k\geq 1$ and any $b\in S\langle a_1,\ldots,a_m\rangle$. Then $L$ is nilpotent.
\end{cor}

\begin{proof}
By \cite[Lemma 5]{MZ}, every element $b\in S\langle a_1,\ldots,a_m\rangle$ is ad-nilpotent of index at most $p^k$. Hence, the claim is an immediate consequence of Theorem \ref{Zelmanov}. 	
\end{proof}

Let $p$ be a prime and $G$ a group. We write $L_p(G)$ for the Lie algebra associated with the Zassenhaus-Jennings-Lazard series $G=D_1\geq D_2\geq\ldots$\, where $$D_i=D_i(G)=\prod_{jp^k\geq i}\gamma_j(G)^{p^k}$$ (see, for instance, \cite[Section 2]{sh00}). Also, we denote by $L(G)$ the subalgebra of $L_p(G)$ generated by $D_1/D_2$.

The next result is a valuable criterion for $L_p(G)$ to satisfy a polynomial identity (see \cite[Theorem 1]{wz92}).

\begin{thm}\label{WZ}
Let $G$ be a group satisfying an identity. Then, for any prime $p$, the Lie algebra $L_p(G)$ satisfies a polynomial identity.
\end{thm}

\section{Proof of Theorem A}

Let $p$ be a prime. A $p$-congruence system in a group $G$ is a descending
chain $G=N_0\supseteq N_1\supseteq N_2\supseteq \ldots$ of normal subgroups of $G$ such that
\begin{itemize}
\item[$(i)$] $G/N_1$ is finite,
\item[$(ii)$] $N_1/N_i$ is a finite $p$-group for all $i>1$, and
\item[$(iii)$] $\bigcap_{i=1}^{\infty}  N_i=\{1\}$.
\end{itemize}

The $p$-congruence system is called of finite rank if there exists a positive integer $r$ such that the rank of $N_1/N_i$ is at most $r$, for all $i$. Recall that the rank of a finite $p$-group $P$ is defined to be the least upper bound of the set $\{d(H)\,|\, H\leq P\}$, where $d(H)$ is the minimal number of generators for $H$. According to a theorem of Lubotzky (see \cite[Theorem B6]{DDMS}), any finitely generated group with a $p$-congruence system of finite rank is linear. 

We say that a group $G$ is residually-$p$ if it is residually a finite $p$-group, that is, for every nontrivial element $x\in G$ there exists a normal subgroup $N$ of $G$ such that $x\notin N$ and $G/N$ is a finite $p$-group.

The following lemma is well-known. For the reader's convenience we provide a proof.

\begin{lem}\label{linear} 
Let $G$ be a finitely generated residually-$p$ group. If $L(G)$ is nilpotent, then $G$ is linear. 
\end{lem}

\begin{proof}
Suppose that $G$ is generated by $m$ elements and $L(G)$ is nilpotent of class $c$. Let $Q$ be a finite homomorphic image of $p$-power order of $G$. Obviously, $Q$ can be generated with $m$ elements and $L(Q)$ is nilpotent of class at most $c$. By \cite[Proposition 1]{KS}, the group $Q$ has a powerful characteristic subgroup of bounded index depending only on $c,m$ and $p$. It follows that the rank of $Q$ is bounded by a function of $c,m$ and $p$ (see \cite[Theorem 2.9]{DDMS}). Since this happens for each finite homomorphic image of $p$-power order of $G$, we conclude that $G$ has a $p$-congruence system of finite rank. Thus $G$ is linear by \cite[Theorem B6]{DDMS}.
\end{proof}

Notice that in a linear group the set of right (or left) Engel elements is contained in the Hirsch-Plotkin radical of $G$ (see \cite[Theorem 0]{Gru}). This leads to an easy but important corollary of Lemma \ref{linear}.

\begin{cor}\label{Gruenberg} 
Let $G$ be a residually-$p$ group generated by finitely many right (or left) Engel elements. If $L(G)$ is nilpotent, then $G$ is nilpotent. 
\end{cor}	

We now deduce the following sufficient condition for nilpotency of a group.

\begin{pro}\label{right}
Let $G$ be a residually-$p$ group such that $L_p(G)$ satisfies a polynomial identity. Suppose that $G$ is generated by finitely many right Engel elements. Then $G$ is nilpotent. 
\end{pro}

\begin{proof}
Let $\{x_1,\dots, x_m\}$ be a finite set of right Engel elements which generate $G$. For any $x_i$, denote by $a_i$ the element $x_iD_2\in L(G)$. Of course $L(G)$ satisfies the same polynomial identity as $L_p(G)$, and it is generated by $a_1,\ldots, a_m$. Take any $b\in S\langle a_1,\ldots, a_m\rangle$ and let $g$ be the group-commutator in $x_i$ having the same system of brackets as $b$. Since each $x_i$ is a right Engel element, there exists $n\geq 1$ such that $[x_i,_n g]=1$ for any $i\in \{1,\ldots,m\}$. Let $p^k$ be the least $p$-power such that $n\leq p^k$. Then in $L(G)$ we have $[a_i,_{p^{k}} b]=0$ and so, by Corollary \ref{MZ}, $L(G)$ is nilpotent. Finally $G$ is nilpotent, by Corollary \ref{Gruenberg}.
\end{proof}

In order to prove Theorem A, we quote a straightforward corollary of \cite[Lemma 2.1]{wi91} (see \cite[Lemma 3.5]{sh00} for a proof).

\begin{lem}\label{wilson}
Let $G$ be a finitely generated residually finite-nilpotent group. For $p$ a prime, denote by $R_p$ the intersection of all normal subgroups of $G$ of finite $p$-power index. If $G/R_p$ is nilpotent for all $p$, then $G$ is nilpotent.
\end{lem}

Recall that, given a residually finite group $G$ satisfying an identity, Theorem A states that the set of right Engel elements of $G$ is contained in the Hirsch-Plotkin radical of $G$.

\begin{proof}[{\bf Proof of Theorem A}]
Let $H$ be the subgroup generated by all right Engel elements of $G$. The claim will follow immediately once it is shown that $H$ is contained in the Hirsch-Plotkin radical of $G$.

Clearly $H$ is normal in $G$. So it is enough to prove that $H$ is locally nilpotent. Take a finitely generated subgroup $K$ of $H$. We may assume that $K$ is generated by finitely many right Engel elements. Since finite groups generated by right Engel elements are nilpotent (see \cite[12.3.7]{Rob}), $K$ is residually finite-nilpotent. Then, by Lemma \ref{wilson}, we may also assume that $K$ is a residually-$p$ group for some prime $p$. Denote by $L_p(K)$ the Lie algebra associated with the Zassenhaus-Jennings-Lazard series of $K$. Now $K$ satisfies the same identity as $G$. Hence, by Theorem \ref{WZ}, $L_p(K)$ satisfies a polynomial identity. Thus $K$ is nilpotent by Proposition \ref{right}, and $H$ is locally nilpotent.
\end{proof}

A group $G$ is called locally graded if every non-trivial finitely generated subgroup of $G$  has a proper subgroup of finite index. Examples of locally graded groups are residually finite groups as well as locally (soluble-by-finite) groups. In \cite{kr} it was shown that any locally graded $n$-Engel group is locally nilpotent. Similarly, one can extend the final part of Theorem A to the class of locally graded groups. 

\begin{cor}
Let $G$ be a locally graded Engel group satisfying an identity. Then $G$ is locally nilpotent.
\end{cor}

\begin{proof}
We may assume that $G$ is finitely generated. Let $R$ be the finite residual of $G$. If $R=1$ then $G$ is residually finite and we are done, by Theorem A. Suppose $R\neq 1$. Of course every section of $G$ satisfies the same identity as $G$. Then $G/R$ is a residually finite group satisfying an identity and therefore nilpotent, by Theorem A. It follows that, for some $d\geq 1$, the $d$-th term $G^{(d)}$ of the derived series of $G$ is a subgroup of $R$. On the other hand, $G/G^{(d)}$ is nilpotent, because it is a soluble group generated by Engel elements \cite[12.3.3]{Rob}. Thus, $R/G^{(d)}$ is a subgroup of the finitely generated nilpotent group $G/G^{(d)}$, so that $R/G^{(d)}$ is also finitely generated. Furthermore, by \cite[Corollary 5]{BSTT13}, $G^{(d)}$ is finitely generated and so $R$ is also finitely generated. Since $G$ is locally graded, there exists a proper subgroup of $R$ with finite index. This implies that the finite residual $S$ of $R$ is a proper subgroup of $R$. By Theorem A, the residually finite group $R/S$ is nilpotent and then $G/S$ is soluble. As above, we conclude that $G/S$ is a finitely generated nilpotent group. Hence $G/S$ is residually finite, which gives $R=S$: a contradiction.
\end{proof}	

\section{Proof of Theorem B}

Let $n\geq 1$ and let $w=w(x_1,\ldots,x_d)$ be an arbitrary word. We denote by $\mathcal W$ the class of all groups $G$ in which the $w$-values are right $n$-Engel and the verbal subgroup $w(G)$ is locally nilpotent. 

Theorem B claims that $\mathcal W$ is a variety. This will be proved applying the following lemma.

\begin{lem}\label{bou} There exists a number $l=l(m,n,w)$ depending only on $m,n,w$ such that if $G\in\mathcal W$ is a nilpotent group generated by $m$ elements which are right $n$-Engel, then the nilpotency class of $G$ is at most $l$.
\end{lem}

\begin{proof}
Suppose that the result is false. Then there exists an infinite sequence $(G_i)_{i\geq 1}$ of nilpotent groups satisfying the hypotheses of the lemma such that the nilpotency class of $G_i$ tends to infinity as $i$ tends to infinity. In particular, each $G_i$ is generated by $m$ elements $x_{i1},\dots,x_{im}$ (not necessarily pairwise distinct) which are right $n$-Engel. Let $C$ be the Cartesian product of the groups $G_i$. For any $i\geq 1$ and $1\leq j\leq m$, let $y_j$ be the element of $C$ whose $i$-th component is equal to $x_{ij}$. Of course each $y_j$ is a right $n$-Engel element in $C$. Put $H=\langle y_1,\ldots,y_m\rangle$. It is easy to verify that $H$ is residually nilpotent and so residually finite, being finitely generated. Moreover, $H$ satisfies the identity $[w(x_1,\ldots,x_d),_n y]\equiv 1$. It follows from Theorem A that $H$ is nilpotent, say of class $c$. Since each $G_i$ is a homomorphic image of $H$, all the subgroups $G_i$ are nilpotent of class at most $c$, which is a contradiction.	
\end{proof}

\begin{proof}[\bf Proof of Theorem B]		
The class $\mathcal W$ is obviously closed with respect to taking subgroups and quotients of its members. Hence, by Birkhoff's theorem (see \cite[2.3.5]{Rob}), it is enough to show that if $G$ is a Cartesian product of groups $G_i$ from $\mathcal W$, then $G\in\mathcal W$. Notice that the $w$-values are right $n$-Engel in $G$. So it remains to prove that the verbal subgroup $w(G)$ is locally nilpotent. 

Let $H$ be a finitely generated subgroup of $w(G)$. Then there exist finitely many $w$-values $w_1,\ldots,w_m$ such that $H\leq\langle w_1,\dots,w_m\rangle$. Assuming $K=\langle w_1,\ldots,w_m \rangle$, let $\pi_i$ be the projection from $K$ to $w(G_i)$. By Lemma \ref{bou} there exists a constant $l$ such that the nilpotency class of $\pi_i(K)$ is at most $l$ for any $i$. Since the intersection of the kernels of $\pi_i$ is trivial, $K$ is nilpotent of class at most $l$. Thus $H$ is nilpotent and $w(G)$ is locally nilpotent.
\end{proof}


\begin{thebibliography}{10}
		
\bibitem{BMTT} R. Bastos, N. Mansuro\u{g}lu, A. Tortora and M. Tota, \textit{Bounded Engel elements in groups satisfying an identity}, Arch. Math. (2018), {https://doi.org/10.1007/s00013-017-1137-x}.

\bibitem{BSTT13} R. Bastos, P. Shumyatsky, A. Tortora and M. Tota, \textit{On groups admitting a word whose values are Engel}, Int. J. Algebra Comput. {\bf 23} no. 1 (2013), 81--89.

\bibitem{DDMS} J.\,D. Dixon, M.\,P.\,F. du Sautoy, A. Mann and D. Segal,
{\it Analytic pro-p-groups}, 2nd edition, Cambridge University Press, Cambridge, 1999.

\bibitem{Gru} K.\,W. Gruenberg, {\em The Engel structure of linear groups},
J. Algebra {\bf 3} (1966), 291--303. 

\bibitem{KS} E.\,I. Khukhro and P. Shumyatsky, {\em Bounding the exponent of a finite group with automorphisms}, J. Algebra {\bf 212} (1999), no. 1, 363--374. 

\bibitem{kr} Y. Kim and A.\,H. Rhemtulla,
{\it On locally graded groups}, in {\it Groups-Korea '94 (Pusan)}, eds. A.\,C. Kim and D.\,L. Johnson (de Gruyter, Berlin, 1995), 189--197.

\bibitem{MZ} C. Mart\'\i nez and E. Zelmanov, {\em On Lie rings of torsion groups}, 
Bull. Math. Sci. {\bf 6} (2016), no. 3, 371--377.
	
\bibitem{Rob} D.\,J.\,S. Robinson,
\textit{A course in the Theory of Groups}, 2nd edition, Springer-Verlag, New York, 1996.

\bibitem{sh00} P. Shumyatsky, {\it Applications of Lie ring methods to group theory}, in {\it Nonassociative algebra and its applications}, eds. R. Costa, A. Grishkov, H. Guzzo Jr. and L.\,A. Peresi, Lecture Notes in Pure and Appl. Math., Vol. 211 (Dekker, New York, 2000), 373--395.
		
\bibitem{STT15} P. Shumyatsky, A. Tortora and M. Tota, \textit{An Engel condition for orderable groups}, Bull. Braz. Math. Soc. (N.S.), {\bf 46} (2015), 461--468.
	
\bibitem{STT16} P. Shumyatsky, A. Tortora and M. Tota, \textit{On locally graded groups with a word whose values are Engel}, Proc. Edinburgh Math. Soc. {\bf 59} (2016), 533--539.
	
\bibitem{STT16_2} P. Shumyatsky, A. Tortora and M. Tota, \textit{On varieties of groups satisfying an Engel type identity}, J. Algebra {\bf 447} (2016), 479--489.

\bibitem{wi91} J.\,S. Wilson, {\it Two-generator conditions for residually finite groups}, Bull. London Math. Soc. {\bf 23} (1991), 239--248.

\bibitem{wz92} J.\,S. Wilson and E.\,I. Zelmanov, {\it Identities for Lie algebras of pro-$p$ groups}, J. Pure Appl. Algebra {\bf 81} (1992), 103--109.
	
\bibitem{ze1} 
E.\,I. Zelmanov, {\it Nil rings and periodic groups}, Korean Math. Soc. Lecture Notes in Math., Seoul, 1992.

\bibitem{ze2} 
E.\,I. Zelmanov, {\it Lie methods in the theory of nilpotent groups}, in {\it Groups '93 Galway/St Andrews}, Vol. 2, 567--585, Cambridge University Press, Cambridge, 1995.
	
\bibitem{ze16} E.\,I. Zelmanov, {\it Lie algebras and torsion groups with identity}, J. Comb. Algebra {\bf 1} (2017), no. 3, 289--340.
	
\end{thebibliography}
\end{document}